\documentclass[11pt]{article}
\usepackage[cp1251]{inputenc} 
\usepackage[dvips]{graphicx}
\usepackage[dvips]{epsfig}
\usepackage{amssymb}
\newtheorem{thm}{Theorem}[section]
 \newtheorem{cor}[thm]{Corollary}
 \newtheorem{lem}[thm]{Lemma}
 \newtheorem{prop}[thm]{Proposition}
 \newtheorem{defn}[thm]{Definition}
 \newtheorem{rem}[thm]{Remark}

\begin{document}

\begin{center}
\Large{{\bf Lifshitz–Kre\u{\i}n trace formula  for Hirsch functional calculus on Banach spaces}}\footnote{This work was financially supported by the  Fund of Fundamental Research of Republic of Belarus. Grant number $\Phi$ 17-082.}
\end{center}

\begin{center}
A. R. Mirotin\\

amirotin@yandex.ru
\end{center}

\textbf{Abstract.}
We give a  simple  definition of a spectral shift function for   pairs of nonpositive operators on Banach spaces and prove  trace formulas of Lifshitz-Kre\u{\i}n type  for a perturbation of an operator monotonic (negative complete Bernstein) function of negative and nonpositive operators on Banach spaces induced by nuclear perturbation of an operator argument. The Lipschitzness  of such functions is also investigated.
The results may be regarded as  a contribution to a perturbation theory for Hirsch functional calculus.

\textbf{Key wards.}
Spectral shift function,  Lifshitz–Kre\u{\i}n trace formula, Hirsch functional calculus, negative operator,   Banach space, perturbation determinant

\section{Introduction}
The trace formula for a trace class perturbation of a self-adjoint operator on Hilbert space was
 proved in a  case of finite-dimensional perturbation by physicist I.~M.~Lifshitz (as a tool for solving some problems in crystal theory)  and in the general case by M.~G.~Kre\u{\i}n \cite{Krein}. Much work has been done during last decades in order to improve and generalize  these results and to get similar formulas (see, e.g., \cite{GMN},  \cite{Pel16}, \cite{AP2}, \cite{MN},  \cite{MNP}, \cite{MNP2}, \cite{Ya} and surveys  \cite{BY}, \cite{AP1} and their references). It should be stressed that all these work deal with  Hilbert spaces only.
The case of operators  on Banach spaces was first considered in \cite{OaM}, \cite{OaMII}, \cite{RM}. Trace formulas of Lifshitz-Kre\u{\i}n type  give an integral representation for the trace of the perturbation   of a function of an operator induced by  a trace class perturbation of an argument using a so called spectral shift function. For applications of such formulas see, e.g., \cite[Ch. 11, 14]{Simon}, \cite[Ch. 8]{Ya}. In this paper, we give a  simple  definition of a spectral shift function for   pairs of nonpositive operators on Banach spaces with nuclear difference and prove  new trace formulas of Lifshitz-Kre\u{\i}n type  for a perturbation of a negative operator monotonic (negative complete Bernstein) function of negative and nonpositive operators on Banach spaces induced by nuclear perturbation of an  argument. The Lipschitzness  of such functions is also investigated.
The results may be regarded as  a contribution to a perturbation theory for Hirsch functional calculus.
\bigskip

\section{ Preliminaries}
In this section we introduce classes of functions and  operators and  briefly  describe  a version of Hirsch functional calculus  we shall use below.

\begin{defn} We say that a function $\varphi$ is  \textit{negative complete Bernstein} and write $\varphi\in \mathcal{OM}_-$ if it is holomorphic in $\mathbb{C}\setminus \mathbb{R}_+,$ satisfies $\mathrm{Im} w \mathrm{Im}\varphi(w)\geq 0$ for $w\in\mathbb{C}\setminus \mathbb{R}_+,$ and such that the limit $\varphi(-0)$ exists
and is real.
\end{defn}

According to \cite[Theorem 6.1]{SSV} this means that $-\varphi(-z)$ is a complete Bernstein function and $\varphi$ has the following integral representation
$$
\varphi(z)=c+bz+\int\limits_{(0,\infty)}\frac{z}{t-z}d\mu(t),\quad z\in \mathbb{C}\setminus (0,+\infty),\eqno(1)
$$
where $c\leq 0, b\geq 0$ and $\mu$ is a unique positive measure that satisfies the condition $\int_{(0,\infty)}d\mu(t)/(1+t)<\infty.$

A lot of  examples of complete Bernstein functions  one can found in \cite{SSV}.

In the sequel unless otherwise stated we  assume for the sake of simplicity  that $c=b=0$ in the integral representation (1) (otherwise one should replace $\varphi(z)$ by $\varphi(z)-c-bz$).

\begin{rem} It is known (see, e.g., \cite[Theorem 12.17]{SSV}), that the families of complete Bernstein and positive operator monotone functions coincide. It follows that the families of negative complete Bernstein and negative operator monotone functions also coincide (we say that a real Borel function $\varphi$ on $(-\infty,0]$ is \textit{negative operator monotone}    if for every bounded self-adjoint operators $A$ and $B$ on a finite or infinite-dimensional real Hilbert space the inequalities $A\leq B\leq O$ imply $\varphi(A)\leq \varphi(B)$). That is why  we denote
the family of negative complete Bernstein functions by $\mathcal{OM}_-.$
\end{rem}

 \begin{defn} We say that (closed, densely defined) operator $A$ on a complex Banach space $X$ is \textit{nonpositive (negative)}  if
 $(0,\infty)$ is contained in $\rho(A),$ the resolvent set of $A$, and
 $$
 M_A:=\sup_{t>0}\|tR(t,A)\|<\infty
 $$
  (respectively   $[0,\infty)\subset \rho(A)$ and
  $$
 M_A:=\sup_{t>0}\|(1+t)R(t,A)\|<\infty)
 $$
  where $R(t,A) = (tI- A)^{-1}$ stands for the resolvent of an operator $A,$ and $Ix=x$ for all $x\in X.$
\end{defn}

So, the operator $A$  is nonpositive (negative)  if and only if $-A$ is non-negative (respectively positive) in a sense of Komatsu \cite{Kom} (see also \cite[Chapter 1]{MS}). We denote by $\mathcal{NP}(X)$ (respectively $\mathcal{N}(X)$) the class of nonpositive (negative) operators on $X.$ (We   deal with  negative  operators instead of positive one because in this form our results are consistent with the multidimensional Bochner-Phillips functional calculus of semigroup generators built in \cite{Mir99} --- \cite{SMZ2011}.) Note  that $A-\varepsilon I\in\mathcal{N}(X)$ for $A\in\mathcal{NP}(X),$  $\varepsilon>0.$

Since every non-negative (closed, densely defined) operator  on  $X$ is sectorial of angle $\omega$  for some $\omega\in (0,\pi)$ (see, e.g.,  \cite[Proposotion 1.2.1]{MS}) every operator  $A\in \mathcal{NP}(X)$ enjoys the following properties ($S_\omega$ denotes the open sector symmetric about the positive real axis with opening angle $2\omega$):

(i) the resolvent set $\rho(A)$ contains the sector $S_\omega,$ $\omega=\arcsin 1/M_A$;

(ii) for every $\omega'\in (0,\omega)$ there is some constant $M'_A\geq M_A$ such that
$$
\|R(\lambda,A)\|\leq \frac{M'_A}{|\lambda|}\quad (\lambda\in S_{\omega'}).
$$

Consequently, if $A$ is  negative then

($\mathrm{i}'$) the resolvent set $\rho(A)$   contains the closure of some set of the form  $S_\theta\cup B_\delta(0)$ ($0<\theta<\arcsin 1/M_A;$ $B_\delta(0)$ stands for the open disc centered at zero of radius $\delta>0$);

($\mathrm{ii}'$)  there is some constant $M'_A\geq M_A$ such that
$$
\|R(\lambda,A)\|\leq \frac{M'_A}{1+|\lambda|}
$$
for $\lambda$ in some neighborhood of the closure of   $S_\theta\cup B_\delta(0).$

  \begin{defn}\cite{BBD} For any function  $\varphi\in \mathcal{OM}_-$ with representing measure $\mu$ and any $A\in \mathcal{NP}(X)$ we put
$$
 \varphi(A)x=\int\limits_{(0,\infty)}AR(t,A)xd\mu(t) \ (x\in D(A)). \eqno(2)
$$
\end{defn}
This operator is closable (see, e.g., \cite{BBD}) and its closure will be denoted by $\varphi(A),$ too.

It is known \cite{Hirsch72} --\cite{Hirsch76} (see also \cite[Theorem 7.4.6]{MS}) that for $\varphi\in \mathcal{OM}_-$ the operator $\varphi(A)$ belongs to $\mathcal{NP}(X)$ ($\mathcal{N}(X)$) if $A\in \mathcal{NP}(X)$ (respectively $A\in\mathcal{N}(X)$).

 \begin{rem} In  the Hirsch functional calculus \cite{Hirsch72} --\cite{Hirsch76} (see also \cite{MS}, \cite{Pys}, \cite{TrudyIM}) functions of the form
$$
f(w)=a+\int\limits_{[0,\infty)}\frac{w}{1+ws}d\lambda(s)
$$
 ($a\geq 0,$  $\lambda$ is a unique positive measure such that $\int_{(0,\infty)}d\lambda(s)/(1+s)<\infty$) are applied to nonnegative operators $T$ on Banach spaces via the formula
 $$
f(T)x=ax+\int\limits_{[0,\infty)}T(I+sT)xd\lambda(s)\quad (x\in D(T)).
$$
Since
$$
f(w)=a+bw+\int\limits_{(0,\infty)}\frac{w}{s^{-1}+w}s^{-1}d\lambda(s),
$$
every such function is complete Bernstein. So, the functional calculus under consideration is in fact a form of Hirsch functional calculus.
\end{rem}
\bigskip

\section{Estimates of perturbations by bounded operators}

This section is devoted to  several  auxiliary results.

\begin{prop} (Cf. \cite{Nab}.) Let $\varphi\in \mathcal{OM}_-$. For any operators $A,B\in \mathcal{NP}(X)$ such that $D(A)\subseteq D(B)$ and
operator $A- B$ is bounded the operator $\varphi(A)-\varphi(B)$ is bounded, too, and the following inequality is valid:
$$
\|\varphi(A)-\varphi(B)\|\leq -(M_A+M_B+M_AM_B)\varphi(-\|A-B\|).
$$
\end{prop}

Proof. Let $A\ne B$. Since $AR(t,A)x=R(t,A)Ax$ for $x\in D(A),$ we have
$$
 (\varphi(A)-\varphi(B))x=\int\limits_{(0,\infty)}(AR(t,A)-BR(t,B))xd\mu(t)\ (x\in D(A)).\eqno(4)
$$
Let $G(t):=AR(t,A)-BR(t,B)\ (t>0).$  The well known equality
$$
AR(t,A)=-I+tR(t,A)\quad (t\in\rho(A))
$$
implies in view of the second resolvent identity that
$$
G(t)=t(R(t,A)-R(t,B))=tR(t,A)(A-B)R(t,B).\eqno(5)
$$
 Hence, by the definition of nonpositive operator, $\|G(t)\|\leq M_AM_B\|A-B\|/t.$ On the other hand, $\|G(t)\|\leq M_A+M_B.$

Now we put in the inequality $(a,b,t>0)$
$$
\min\left\{a,\frac{1}{t}\right\}\leq\frac{1+ab}{t+b}
$$
$a=(M_A+M_B)/(M_AM_B\|A-B\|), b=\|A-B\|,$ and obtain
$$
\|G(t)\|\leq \frac{(M_A+M_B+M_AM_B)\|A-B\|}{t+\|A-B\|}.\eqno(6)
$$
It follows that the Bochner integral $\int_{(0,\infty)}G(t)d\mu(t)$ exists with respect to the operator norm, the operator  $\varphi(A)-\varphi(B)$ is bounded, and by  formulas (4) and (5)
$$
 \varphi(A)-\varphi(B)=
  \int\limits_{(0,\infty)}(R(t,A)-R(t,B))td\mu(t).\eqno(7)
$$

 Moreover (6) yields
$$
\|\varphi(A)-\varphi(B)\|\leq(M_A+M_B+M_AM_B)\|A-B\|\int\limits_{(0,\infty)}\frac{d\mu(t)}{t+\|A-B\|}=
$$
$$
-(M_A+M_B+M_AM_B)\varphi(-\|A-B\|).
$$

 \begin{cor} For the function  $\varphi\in \mathcal{OM}_-$ (with $b=c=0$) the following statements are equivalent:

1)  $\varphi$ is operator Lipschitz in any class $\mathcal{NP}_c(X):=\{A\in \mathcal{NP}(X):M_A\leq c\},$ $c=\mathrm{const}$
(i.e. $\|\varphi(A)-\varphi(B)\|\leq L\|A-B\|$ for all $A,B\in \mathcal{NP}_c(X)$ such that $A-B$ is bounded);

2) $\varphi'(-0)\ne \infty;$

3) $\int_{(0,\infty)}d\mu(t)/t<\infty;$

4) $|\varphi(x)|\leq L_1|x|$  for all $x\in (-\infty,0];$

5)  $\varphi$ is Lipschitz on $(-\infty,0]$.
\end{cor}

Proof. The equivalence 1)$\Leftrightarrow$ 4) follows from
proposition  3.1.
Since by the Monotone Convergence Theorem $\varphi'(-0)=\int\limits_{(0,\infty)}d\mu(t)/dt,$ we get the statement 2)$\Leftrightarrow$ 3).
The implications 4) $\Rightarrow$ 2), 1) $\Rightarrow$ 5), and 5) $\Rightarrow$ 4) are obvious.  And finally, 3) $\Rightarrow$ 4), since for all $x\leq 0$
$$
|\varphi(x)|\leq \int\limits_{(0,\infty)}\frac{|x|t}{t-x}\frac{d\mu(t)}{t}\leq \left(\int\limits_{(0,\infty)}\frac{d\mu(t)}{t}\right)|x|.
$$

The functional calculus under consideration satisfies the following \textit{stability property}.

 \begin{cor} Let $\varphi\in \mathcal{OM}_-.$  For any sequences of operators $A_n,B_n\in \mathcal{NP}(X)$ such that $D(A_n)\subseteq D(B_n),$  $M_{A_n}, M_{B_n}<\mathrm{const}$ and
 $\|A_n- B_n\|\to 0$ we have $\|\varphi(A_n)-\varphi(B_n)\|\to 0$ ($n\to\infty$).
\end{cor}

\begin{prop} If under assumptions of proposition  3.1 the operators $A-B$ and $R(t,B)$ commute, then
for any $x \in D(A), \|x\|=1$ the following inequality is valid:
$$
\|(\varphi(A)-\varphi(B))x\|\leq -(M_A+M_B+M_AM_B)\varphi(-\|(A-B)x\|).
$$
\end{prop}

This proposition can be proved in just the same way as proposition  3.1.

 \begin{cor} Let $A\in \mathcal{NP}(X),$ $\varphi\in \mathcal{OM}_-.$ Then for any $x \in D(A)$ the following inequalities are valid:

1) (cf. \cite[Corollary 13.8]{SSV})
$$
\|\varphi(A)x\|\leq -(2M_A+1)\varphi(-\|Ax\|)\ (\|x\|=1);
$$

2)
$$
\|\varphi(A)x\|\leq (2M_A+1)\varphi'(-0)\|Ax\|\ (x\ne 0).
$$
\end{cor}

Proof. 1) It is a special case of proposition  3.4 for $B = O.$

2) Note that the function $\varphi(-s)/(-s)$ decreases on $\{s>0\}$ in view of formula (1). It  follows that  $\varphi'(-0)\geq \varphi(-s)/(-s)$ for all $s>0$ (we assume that  $c=\varphi(0)=0$).
In particular, $-\varphi(-\|Ax\|)\leq \varphi'(-0)\|Ax\|)$ and therefore 1) implies 2).

In what follows $(\mathcal{I},\|\cdot\|_{\mathcal{I}})$ stands for an operator ideal in $X.$ That means that $\mathcal{I}$ is a two-sided ideal of the algebra $\mathcal{L}(X)$ of bounded operators on  $X,$ the ideal $\mathcal{I}$ is complete with respect to  the norm $\|\cdot\|_{\mathcal{I}},$ and  the following conditions hold: $\|ASB\|_{\mathcal{I}}\leq \|A\|\|S\|_{\mathcal{I}}\|B\|,$ $\|S\|\leq \|S\|_{\mathcal{I}}$
 for all  $A, B\in \mathcal{L}(X)$ and $S\in \mathcal{I}$ (the case $\mathcal{I} = \mathcal{L}(X)$ is not excluded, and is of interest).

\begin{prop} (Cf. \cite{Nab}.) Let $\varphi\in \mathcal{OM}_-,$ $\varphi'(-0)\ne \infty.$  For any operators $A,B\in \mathcal{NP}(X)$ such that $D(A) \subseteq D(B)$ and
$A - B$ belongs to $\mathcal{I},$ the operator $\varphi(A)-\varphi(B)$ also belongs to $\mathcal{I}$ and satisfies the inequality
$$
\|\varphi(A)-\varphi(B)\|_\mathcal{I}\leq M_AM_B\varphi'(-0)\|A-B\|_\mathcal{I}.
$$
\end{prop}

Proof. Formula (5) shows that $G(t)\in \mathcal{I}$  and
$$
\|G(t)\|_{ \mathcal{I}}\leq \frac{M_AM_B\|A-B\|_{ \mathcal{I}}}{t}
$$
for $t>0.$ Since $\int_{(0,\infty)}d\mu(t)/dt=\varphi'(-0),$ it follows that the Bochner integral in (7) exists with respect to the norm $\|\cdot\|_{\mathcal{I}}$ and
$$
\|\varphi(A)-\varphi(B)\|_\mathcal{I}\leq M_AM_B\|A-B\|_\mathcal{I}\int\limits_{(0,\infty)}\frac{d\mu(t)}{t} = M_AM_B\varphi'(-0)\|A-B\|_\mathcal{I}.
$$

\begin{cor} Let $\varphi\in\mathcal{ OM}_-,$ $\varphi'(-0)\ne \infty$. For any sequences of operators $A_n,B_n\in \mathcal{NP}(X)$ such that $D(A_n)\subseteq D(B_n),$  $M_{A_n}, M_{B_n}<\mathrm{const},$ $A_n- B_n\in \mathcal{I},$ and
 $\|A_n- B_n\|_\mathcal{I}\to 0$ we have $\|\varphi(A_n)-\varphi(B_n)\|_\mathcal{I}\to 0$ ($n\to\infty$).
\end{cor}

\begin{cor} Let $U$ be an automorphism of the space $X,$  $\varphi\in \mathcal{OM}_-,$ $\varphi'(-0)\ne \infty,$ $A\in \mathcal{NP}(X).$ If $[A,U] \in\mathcal{ I},$ then  $[\varphi(A),U]\in \mathcal{I}$ and
$$
\|[\varphi(A),U]\|_\mathcal{I}\leq M_A^2\varphi'(-0)\|[A,U]\|_\mathcal{I}.
$$
\end{cor}

Proof. Note that $R(t,UAU^{-1}) = UR(t,A)U^{-1}$ $(t>0).$ Then $UAU^{-1}\in \mathcal{NP}(X),$  $M_{UAU^{-1}} = M_A,$ and $\varphi(UAU^{-1})x = U\varphi(A)U^{-1}x$  for all $x\in D(A).$ Since $D(A)$ is a core for (closed) operators standing in both sides of the last equality,  $\varphi(UAU^{-1}) = U\varphi(A)U^{-1}.$  It follows that  $[\varphi(A), U] = (\varphi(A)-\varphi(UAU^{-1}))U.$ Since
$A- UAU^{-1} = [A,U]U^{-1}\in \mathcal{I},$  proposition  3.6 yields
$$
\|[\varphi(A),U]\|_\mathcal{I}\leq  \|\varphi(A)-\varphi(UAU^{-1})\|_\mathcal{I}\leq
$$
$$
 M_A^2\varphi'(-0)\|A-UAU^{-1}\|_\mathcal{I}\leq
   M_A^2\varphi'(-0)\|[A,U]\|_\mathcal{I}.
$$
\bigskip

\section{Lifshitz-Kre\u{\i}n trace formula}

\subsection{Main results}
In this subsection we introduce a spectral shift function and prove an analog of Lifshitz-Kre\u{\i}n trace formula for  pairs of negative and nonpositive operators on a Banach space.

First note that the function
$\psi_\lambda(s):=\log \lambda-\log(\lambda-s)$ ($\lambda>0$) belongs to $\mathcal{OM}_-$ \cite[Example 3]{MirSF}, \cite{SSV}. So, for $A\in \mathcal{NP}(X), \lambda>0$ we can put
$$
\log(\lambda I-A):=(\log\lambda) I-\psi_\lambda(A).
$$
Note also that for $A, B\in \mathcal{NP}(X)$ such that $A-B$  is nuclear and $\lambda>0$ the operator
$$
\log(\lambda I-A)-\log(\lambda I-B)=\psi_\lambda(B)-\psi_\lambda(A)
$$
is nuclear by proposition  3.6.
(Recall that operator on $X$ is \textit{ nuclear} if it is representable as the sum of absolutely convergent in operator norm series  of rank one operators, see, e.g., \cite[p. 64]{DF}.)

\begin{defn} Let the Banach space $X$ has the approximation property (see, e.g., \cite{DF}). For $A, B\in \mathcal{NP}(X), \lambda>0$ such that $D(A)\subseteq D(B)$ and $A-B$  is nuclear define
\textit{the spectral shift function for the  pair} $(A,B)$ for $\lambda>0$ by
$$
\xi_{A,B}(\lambda)=\mathrm{tr}(\log(\lambda I-A)-\log(\lambda I-B)).
$$
\end{defn}

\begin{thm}  Let the Banach space $X$ has the approximation property.  Let  $A$ and $B$
be negative operators on $X$ such that $D(A)\subseteq D(B)$ and $A-B$ is nuclear.
There exists an analytic  continuation of the  spectral shift  function $\xi_{A,B}$ into the  closure of
some domain $\Omega_{A,B}$ of the form  $S_\theta\cup B_\delta(0)$ ($\theta\in (0,\pi/2)$)
 such that for every negative operator monotone function $\varphi,$ $\varphi'(-0)\ne \infty,$ with the property
$$
\int\limits_0^{\infty}\frac{|\varphi(-x)|}{1+x^2}dx<\infty\eqno(\ast)
$$
the following trace formula holds:
$$
\mathrm{tr}(\varphi(A)-\varphi(B))=\frac{1}{2\pi i}\int\limits_{\Gamma_{A,B}}\xi_{A,B}(z)\varphi'(z)dz\eqno(LK)
$$
where $\Gamma_{A,B}$ denotes the positive oriented  boundary of $\Omega_{A,B}.$

Conversely, if the formula $(LK)$ holds  for every pare $(A,B)$ of negative operators
 on the one-dimensional complex space,  the function $\varphi$ satisfies  the condition $(\ast).$
\end{thm}

Proof. Proposition  3.6 implies that the operator $\varphi(A)-\varphi(B)$ belongs to the ideal $\mathfrak{S}_1=\mathfrak{S}_1(X)$ of nuclear operators on $X.$ Moreover, since by the second resolvent identity and condition $\mathrm{(ii')}$
$$
\|R(\zeta,A)-R(\zeta,B)\|_{\mathfrak{S}_1}\leq \frac{M'_AM'_B\|A-B\|_{\mathfrak{S}_1}}{1+|\zeta|^2}, \quad (\zeta\in\overline{S_\theta}, \theta\in(0,\pi/2)) \eqno(8)
$$
the Bochner integral in (7) converges with respect to the nuclear norm and
$$
\mathrm{tr}(\varphi(A)-\varphi(B))=\int\limits_0^\infty \eta_{A,B}(t)td\mu(t),\eqno(9)
$$
where the function
$$
\eta_{A,B}(z):=\mathrm{tr}(R(z,A)-R(z,B))
$$
 is holomorphic in $\rho(A)\cap\rho(B).$ Indeed, fix $z_0\in \rho(A)\cap\rho(B).$ For some neighborhood of  $z_0$ we have
$$
R(z,A)=\sum\limits_{n=0}^\infty (z-z_0)^nA_n, R(z,B)=\sum\limits_{m=0}^\infty (z-z_0)^mB_m
$$
(both series with operator coefficients converge absolutely in the operator norm). Hence
$$
R(z,A)(A-B)R(z,B)=\sum_{n=0}^\infty\sum\limits_{m=0}^\infty (z-z_0)^{n+m}A_n(A-B)B_m,
$$
where the series in the right-hand side converges in the nuclear norm due to the inequality
$$
\|(z-z_0)^{n+m}A_n(A-B)B_m\|_{\mathfrak{S}_1}\leq |z-z_0|^{n+m}\|A_n\|\|B_m\|\|(A-B)\|_{\mathfrak{S}_1}.
$$

The  set $\rho(A)\cap\rho(B)$ contains the closure of some set $\Omega_{A,B}$ of the form  $S_\theta\cup B_\delta(0)$ ($0<\theta<\arcsin1/M_A, \arcsin1/M_B$) such that the condition ($\mathrm{ii}'$) from the Preliminaries holds.

Note, that for $s<0$
$$
\psi_\lambda(s)=\int\limits_\lambda^\infty\frac{s}{t-s}\frac{dt}{t}.
$$
So, by formula (9) ($\lambda>0$),
$$
\xi_{A,B}(\lambda)=\mathrm{tr}(\psi_\lambda(B)-\psi_\lambda(A))=-\int\limits_\lambda^\infty\eta_{A,B}(t)dt.
$$
Let $L_z$  denotes the ray in $ \overline{S_\theta}$ that starts at $z\in \overline{S_\theta}$ and has a slope $\tan\theta.$ Since  $|\eta_{A,B}(z)|\leq \mathrm{C}/(1+|z|^2)$ ($C=M'_AM'_B\|A-B\|_{\mathfrak{S}_1}$) by  (8) and  $\eta_{A,B}$ is holomorphic in $\rho(A)\cap\rho(B),$ we have for $\lambda>0$ by the Cauchy Theorem
$$
\int\limits_\lambda^\infty\eta_{A,B}(t)dt=\int\limits_{L_\lambda}\eta_{A,B}(\zeta)d\zeta.
$$
Thus the formula
$$
\xi_{A,B}(z)=-\int\limits_{L_z}\eta_{A,B}(\zeta)d\zeta \quad (z\in \overline{S_\theta}) \eqno(10)
$$
 gives   the analytic continuation of  $\xi_{A,B}$   into  the closure $\overline{S_\theta}$ of $S_\theta$ such that $\xi'_{A,B}(z)=\eta_{A,B}(z).$ We claim  that the integral in (10) converges and  the following estimate holds
$$
|\xi_{A,B}(z)|\leq\frac{C}{\mathrm{Re}z}\quad(z\in\overline{\Omega_{A,B}}, \mathrm{Re}z>0).\eqno(11)
$$
Indeed, let $\zeta\in L_z,$ $x:=\mathrm{Re}\zeta.$ As shown above
$|\eta_{A,B}(\zeta)|\leq C/(1+|\zeta|)^2\leq C/x^2.$ So, in view of $|d\zeta|=dx/\cos\theta$ we have
$$
|\xi_{A,B}(z)|\leq \int\limits_{L_z}|\eta_{A,B}(\zeta)||d\zeta| \leq C\int\limits_{\mathrm{Re}z}^\infty\frac{dx}{x^2} =\frac{C}{\mathrm{Re}z}.
$$

We denote also by $\xi_{A,B}$ the antiderivative for $\eta_{A,B}$ which is the analytic continuation of $\xi_{A,B}$ from $\overline{S_\theta}$ to some neighborhood of the closure of $\Omega_{A,B}.$
  Then for every $t\geq 0$ and $E>t$ the Cauchy formula holds:
 $$
\eta_{A,B}(t)=\frac{1}{2\pi i}\int\limits_{\partial G_E}\frac{\xi_{A,B}(z)}{(t-z)^2}dz, \eqno(12)
$$
where $G_E=\{z\in \Omega_{A,B}:\mathrm{Re }z\leq E\}$ and $\partial G_E$ denotes the positive oriented boundary of $G_E.$

Consider the segment  $T_E:=\{z\in \Omega_{A,B}:\mathrm{Re }z= E\}.$ Then
 $$
\lim\limits_{E\to\infty}\int\limits_{T_E}\frac{\xi_{A,B}(z)}{(t-z)^2}dz=0.\eqno(13)
$$
Indeed, taking into account that the length of the segment $T_E$ is $2E\tan\theta,$ we get in view of (11)
 $$
\left|\int\limits_{T_E}\frac{\xi_{A,B}(z)}{(t-z)^2}dz\right|\leq \frac{2\tan\theta}{(E-t)^2},
$$
and (13) follows. Putting (12) and (13) together we obtain for $t\geq 0$
 $$
\eta_{A,B}(t)=\frac{1}{2\pi i}\int\limits_{\Gamma_{A,B}}\frac{\xi_{A,B}(z)}{(t-z)^2}dz, \eqno(14)
$$
where $\Gamma_{A,B}$ denotes the positive oriented boundary of $\Omega_{A,B}.$ In turn, putting together (9) and (14), we get in view of the Fubini Theorem that
$$
\mathrm{tr}(\varphi(A)-\varphi(B))=\int\limits_0^\infty \frac{1}{2\pi i}\int\limits_{\Gamma_{A,B}}\frac{\xi_{A,B}(z)}{(t-z)^2}dz td\mu(t)=
$$
$$
\frac{1}{2\pi i}\int\limits_{\Gamma_{A,B}}\xi_{A,B}(z)\int\limits_0^\infty \frac{td\mu(t)}{(t-z)^2}dt=
\frac{1}{2\pi i}\int\limits_{\Gamma_{A,B}}\xi_{A,B}(z)\varphi'(z)dz.
$$
To complete  the proof of (LK) it remains to legitimate the application of   Fubini Theorem. To this end we are going to deduce from $(\ast)$ the convergence of integrals
$$
I_k:=\int\limits_{\Gamma_{k}}\int\limits_0^\infty \frac{td\mu(t)}{|t-z|^2}dt|\xi_{A,B}(z)||dz|\quad (k=0,1,2)
$$
where $\Gamma_0:=\{z\in \Gamma_{A,B}:\mathrm{Re}z\leq h\}=\partial \Omega_{A,B}\cap\partial B_\delta(0)$ is the arc of the sircle $\partial B_\delta(0),$ and  $\Gamma_{1,2}:=\{z\in \Gamma_{A,B}: \arg z=\pm\theta\}.$

First of all note that for all  $t\in \mathbb{R}_+$ and $z\in \mathbb{C}$ such that $|\arg z|\geq \theta$
$$
|z-t|^2=|z|^2+t^2-2t|z|\cos(\arg z)\geq (|z|^2+t^2)(1-\cos\theta). \eqno(15)
$$
Therefore  ($x:=\mathrm{Re}z$)
$$
\int\limits_0^\infty \frac{td\mu(t)}{|t-z|^2}dt\leq \frac{1}{1-\cos\theta}\int\limits_0^\infty \frac{td\mu(t)}{t^2+|z|^2}dt\leq \frac{1}{1-\cos\theta}\int\limits_0^\infty \frac{td\mu(t)}{t^2+x^2}dt
$$
 ($|\arg z|\geq \theta$  for $z\in \Gamma_{A,B}$).
Since for $z\in \Gamma_{1,2}$ we have $x\geq h$ for some constant $h>0$ and $|dz|=dx/\cos\theta,$ it follows by virtue of formula (11), that
$$
I_{1,2}\leq \mathrm{const}\int\limits_h^\infty \frac{1}{x}\int\limits_0^\infty \frac{td\mu(t)}{t^2+x^2}dt.
$$
Moreover,
$$
\int\limits_1^\infty \frac{1}{x}\int\limits_0^\infty \frac{td\mu(t)}{t^2+x^2}dt=\int\limits_0^\infty\int\limits_1^\infty \frac{dx}{x(x^2+t^2)}td\mu(t)=\frac{1}{2}\int\limits_0^\infty\frac{\log(1+t^2)}{t}d\mu(t),
$$
and the condition $(\ast)$ implies that the last integral converges, because
$$
\int\limits_0^{\infty}\frac{|\varphi(-x)|}{1+x^2}dx=\int\limits_0^{\infty}\left(\int\limits_0^{\infty}\frac{xd\mu(t)}{t+x} \right)\frac{dx}{1+x^2}=
$$
$$
\int\limits_0^{\infty}\int\limits_0^{\infty}\frac{x}{(1+x^2)(x+t)}dxd\mu(t)=\frac{1}{2}\int\limits_0^{\infty}\frac{2t\log t+\pi}{1+t^2}d\mu(t).
$$

To prove the convergence of $I_0$, note that for $z\in \Gamma_0$ formula (15) yields $|z-t|^2\geq (h^2+t^2)(1-\cos\theta),$ and hence
$$
\int\limits_0^{\infty}\frac{td\mu(t)}{|z-t|^2}\leq\frac{1}{1-\cos\theta}\int\limits_0^{\infty}\frac{td\mu(t)}{h^2+t^2}<\infty.
$$
Since $\xi_{A,B}$ is bounded on  $\Gamma_0,$ it follows that $I_0<\infty.$  This completes the proof of (LK).

To prove the last statement of the theorem, assume that (LK) holds for operators $A=-I, B=-2I$ on the one-dimensional complex space. Since by (10) $\xi_{A,B}(z)=\log\frac{z+1}{z+2}$ ($\log 1=0$), formula (LK) implies the convergence of the integral
$$
\int\limits_{\Gamma_1}\log\frac{z+1}{z+2}\varphi'(z)dz=
\left.\log\frac{z+1}{z+2}\varphi(z)\right|_{z_0}^\infty-\int\limits_{\Gamma_1}\frac{\varphi(z)}{(z+1)(z+2)}dz \eqno(16)
$$
($z_0$ denotes the origin of $\Gamma_1$). Here integration by parts is legal since (see \cite[p. 76]{SSV})
$\lim\limits_{r\to \infty}\varphi(re^{i\theta})/r=0$
and therefore
$$
\lim\limits_{z\to \infty, z\in \Gamma_1}\log\frac{z+1}{z+2}\varphi(z)=\lim\limits_{z\to \infty, z\in \Gamma_1}(-\frac{1}{z+2}-\frac{1}{2(z+2)^2}-\dots)\varphi(z)=0.
$$
Now it follows from (16) that the integral $\int\limits_{\Gamma_1}\varphi(z)/(z^2+1)dz$ converges as well. For $R>0$ consider the curves
$$
\Gamma_{1, R}:=\{z\in \Gamma_1: |z|\leq R\}=\{z: \delta\leq |z|\leq R, \arg z=\theta\},
$$
$$
C_R:=\{z: |z|=R, \theta\leq\arg z\leq \pi\},\ \gamma_\delta:=\{z: |z|=\delta, \theta\leq\arg z\leq \pi\}.
$$
 By the Cauchy Theorem
$$
\int\limits_{\Gamma_{1, R}}+\int\limits_{C_R}+\int\limits_{[-R,-\delta]}+\int\limits_{\gamma_\delta}\frac{\varphi(z)}{z^2+1}dz=0.
$$
So, it remains to prove that
$$
\lim_{R\to\infty}\int\limits_{C_R}\frac{\varphi(z)}{z^2+1}dz=0
$$
or, equivalently,
$$
\lim_{R\to\infty}\int\limits_{C_R}\left(\int\limits_0^\infty\frac{zd\mu(t)}{z-t}\right)\frac{dz}{z^2+1}=0.\eqno(17)
$$
For the proof we consider the following integral
$$
\int\limits_0^\infty\left(\int\limits_{C_R}\frac{|z|}{|z^2+1||z-t|}|dz|\right)d\mu(t).\eqno(18)
$$
 For $z\in C_R$ formula (15) yields $|z-t|^2\geq (R^2+t^2)(1-\cos\theta),$ and therefore for $R\geq 2$ we have
$$
\int\limits_{C_R}\frac{|z|}{|z-t||z^2+1|}|dz|\leq \frac{1}{\sqrt{1-\cos\theta}}\frac{R}{\sqrt{R^2+t^2}(R^2-1)}\pi R\leq \frac{1}{\sqrt{1-\cos\theta}}\frac{4\pi}{R+t}.
$$
Since, by the Monotone Convergence Theorem
$$
\lim_{R\to\infty}\int\limits_0^\infty\frac{d\mu(t)}{R+t}=0,
$$
the formula (17) follows.

\begin{rem}
 The function
$$\varphi(s)=\frac{s\log(-s)-s-1}{\log^2(-s)}\quad (s<0)
$$
belongs to $\mathcal{OM}_-$ \cite[p. 337]{SSV},  satisfies   $c=b=0,$ and $\varphi'(-0)=0$ but the condition $(\ast)$ does not fulfilled for this function. So,  the formula (LK) does not hold for $\varphi.$
\end{rem}

\begin{rem}
 It was shown in \cite{Pel16}
that Lifshitz-Kre\u{\i}n trace formula  holds for arbitrary pairs  of not necessarily bounded
self-adjoint operators with trace class difference if and only if the corresponding function is operator Lipschitz (in the class of self-adjoint operators).  The condition $\varphi'(-0)\ne \infty$ guarantee (see corollary 3.2) that the function $\varphi$ is operator Lipschitz in any class $\mathcal{NP}_c(X)$. So,  the result of theorem 4.2 is consistent with the result of V.~V.~ Peller mentioned above (nonpositive operators on Hilbert space $H$ belong to $\mathcal{NP}_1(H)$).
\end{rem}

\begin{cor} (Cf. \cite[formula (2.10)]{Krein}.) \textit{Under the conditions of theorem 4.2 we have for} $\lambda\geq 0$
$$
\mathrm{tr}(R(\lambda,A)-R(\lambda,B))=\frac{1}{2\pi i}\int\limits_{\Gamma_{A,B}}\frac{\xi_{A,B}(z)}{(\lambda-z)^2}dz.
$$
\end{cor}

Proof.
 This follows from formula (14).

Now we generalize the notion of a perturbation determinant to the case of operators on Banach spaces (c.f. \cite{OaMII}). First note that for $A, B\in \mathcal{NP}(X)$ such that $D(A)\subseteq D(B)$ and $A-B$  is nuclear  $\xi_{A,B}$ has an analytic continuation to some sector $S_\theta$ with $\theta\in(0,\pi/2)$ (see theorem 4.12 below).

\begin{defn} For $A, B\in \mathcal{NP}(X)$ such that $D(A)\subseteq D(B)$ and $A-B$  is nuclear   define the
\textit{ perturbation determinant  for the pair}  $(A,B)$ for $z\in S_\theta$ as follows:
$$
\Delta_{B/A}(z)=\exp(\xi_{A,B}(z))
$$
(for negative operators one can take $z\in \Omega_{A,B}$).
\end{defn}

 Since the formula (11) remains true  for nonpositive operators, we have
$$
|\Delta_{B/A}(z)-1|\leq \exp|\xi_{A,B}(z)|-1\leq \exp\left(\frac{C}{\mathrm{Re}z}\right)-1.
$$
So, $\Delta_{B/A}(z)$ belongs to the open  right half-plane for all $z\in S_\theta$ with sufficiently large $\mathrm{Re}z.$ Since (again by the formula (11)) $\xi_{A,B}(t)\to 0$ as $t\to +\infty,$ it follows that for such $z$

$$
\xi_{A,B}(z)=\log\Delta_{B/A}(z)
$$
where $\log$ stands for the branch of the logarithm in the right half-plane that satisfies $\log 1=0.$

\begin{cor} (Cf. \cite[formula (3.8)]{Krein}.) \textit{Under the conditions of theorem 4.2 we have for sufficiently large} $\lambda> 0$
$$
\xi_{A,B}(\lambda)=\log\Delta_{B/A}(\lambda)=\frac{1}{2\pi i}\int\limits_{\Gamma_{A,B}}\frac{\xi_{A,B}(z)}{z-\lambda}dz.
$$
\end{cor}

Proof.
 This follows from (LK) with $\varphi=\psi_\lambda.$

 \begin{cor} (Cf. \cite[formula (3.17)]{Krein}.) \textit{Under the conditions of theorem 4.2 we have for all $z\in \Omega_{A,B}$}
$$
\frac{\Delta'_{B/A}(z)}{\Delta_{B/A}(z)}=\mathrm{tr}(R(z,A)-R(z,B)).
$$
\end{cor}

Proof.
Indeed, $\Delta'_{B/A}(z)/\Delta_{B/A}(z)=\xi'_{A,B}(z)=\eta_{A,B}(z)$ (formula (10)).

The following  lemma will be  useful.

\begin{lem} \textit{Let $A \in \mathcal{N}(X).$ Then  $A+V \in \mathcal{N}(X)$ and $\rho(A+V)\supset \rho(A)$ for every $V\in \mathcal{L}(X)$ such that $\|V\|<1/M'_A.$ In this case one can take $M'_{A+V}= M'_A(1-M'_A\|V\|)^{-1}.$}
\end{lem}

Proof. First note that  $\|V\|<1/M'_A\leq \|R(\lambda,A)\|^{-1}$ for all  $\lambda$ in some neighborhood $\mathcal{O}$ of the closure of   $S_\theta\cup B_\delta(0)$  due to ($\mathrm{ii}'$). It follows in view of \cite[Remark IV.3.2]{Kato}   that $\rho(A+V)\supset \rho(A)\supset\mathbb{ R}_+.$

Applying \cite[Theorem IV.1.16, Remark IV.1.17]{Kato} we have for  $\lambda\in \mathcal{O}$
$$
\|R(\lambda,A+V)-R(\lambda,A)\|\leq \frac{\|V\|\|R(\lambda,A)\|^2}{1-\|V\|\|R(\lambda,A)\|},\eqno(19)
$$
since $\|V\|\|R(\lambda,A)\|< \|R(\lambda,A)\|/M'_A\leq 1.$

Thus, using the condition ($\mathrm{ii}'$) once more  we obtain for  $\lambda\in \mathcal{O}$
$$
\|R(\lambda,A+V)\|\leq \|R(\lambda,A)\|+\frac{\|V\|\|R(\lambda,A)\|^2}{1-\|V\|\|R(\lambda,A)\|}=
$$
$$
\frac{\|R(\lambda,A)\|}{1-\|V\|\|R(\lambda,A)\|}\leq
\frac{\|R(\lambda,A)\|}{1-\|V\|M'_A}\leq \frac{M'_A(1-M'_A\|V\|)^{-1}}{1+|\lambda|}.
$$

Now we are in position to prove a formula for the spectral shift function.
In the following two theorems we assume that the nuclear operator $A-B$ has the
 form
$$
A-B=\sum\limits_{j=1}^\infty\ell_j\otimes v_j\eqno(20)
$$
where $\ell_j \otimes v_j(x):=\ell_j(x)v_j\ (\ell_j\in X'; v_j,x\in X),$ the tensor product of a linea functional $\ell_j$ and vector  $v_j,$ and $\sum_{j=1}^\infty\|\ell_j\|\|v_j\|<\infty.$ In this case $\mathrm{tr}(A-B):=\sum_{j=1}^\infty\ell_j(v_j).$

\begin{thm} (Cf., e.g., \cite[(3.11) and (3.4)]{BY}.) Let the Banach space $X$ has the approximation property. For any operators $A, B\in \mathcal{N}(X)$ such that $D(A)\subseteq D(B)$ and the
operator $A - B$ has the form (20) the following equality is valid for $\lambda\in \Omega_{A,B}$  with sufficiently large $|\lambda|$:
$$
\xi_{A,B}(\lambda)=\sum\limits_{k=1}^\infty\log(1-\ell_k(R(\lambda,A_{k-1})v_k)),
$$
where $A_0:=B,$ $A_k:=B+S_k,$ $S_k:=\sum_{j=1}^k\ell_j\otimes v_j$ ($k\in \mathbb{N}$),  and $\log$ denotes the branch of the logarithm in the right half-plane that satisfies $\log 1=0;$ the series converges absolutely.
\end{thm}

Proof. Let $R_n:=\sum_{j=n+1}^\infty\ell_j\otimes v_j.$ Choose such $N\in \mathbb{N}$  that $\|R_n\|<1/M'_A$ for all $n\geq N.$ Then the operator $A_n=A-R_n$ belongs to $\mathcal{N}(X)$ for all $n\geq N$  by lemma 4.9.

Let $\lambda_0:=M'_B\sum_{j=1}^\infty\|\ell_j\|\|v_j\|.$  First note that $R(\lambda,A_k)=((\lambda I-B)-S_k)^{-1}$ exists for $\lambda\in \Omega_{A,B},$ $|\lambda|>\lambda_0$ because  $\|S_k\|\leq \sum_{j=1}^\infty\|\ell_j\|\|v_j\| <1/\|R(\lambda, B)\|$  for such $\lambda$ due to ($\mathrm{ii}'$).

We claim  that for all $n\geq N$
$$
\Delta_{B/A_n}(\lambda)=\prod\limits_{k=1}^n(1-\ell_k(R(\lambda,A_{k-1})v_k))\quad (\lambda\in\rho(B), |\lambda|>\lambda_0).\eqno(21)
$$

First of all, using an approach by Kre\u{\i}n we compute $R(\lambda,A_1)$  for $\lambda\in\rho(B).$ In this case the equation
$$
\lambda x-A_1x=y\quad (y\in X, \lambda\in\rho(B),)
$$
has the form
$$
(\lambda-B)x-\ell_1(x)v_1=y,
$$
or, equivalently,
$$
x=\ell_1(x)R(\lambda,B)v_1+R(\lambda,B)y.
$$
If we denote
$$
a=\ell_1(x) \eqno(22)
$$
then
$$
x=aR(\lambda,B)v_1+R(\lambda,B)y.
$$
Substituting this into (22) we get
$$
a=\frac{\ell_1(R(\lambda,B)y)}{1-\ell_1(R(\lambda,B)v_1)}.
$$
It follows that
$$
R(\lambda,A_1)y=x=\frac{\ell_1(R(\lambda,B)y)}{1-\ell_1(R(\lambda,B)v_1)}R(\lambda,B)v_1+R(\lambda,B)y.
$$
So,
$$
R(\lambda,A_1)-R(\lambda,B)=\ell\otimes (R(\lambda,B)v_1),
$$
the one-dimensional operator, where
$$
\ell(x):=\frac{\ell_1(R(\lambda,B)x)}{1-\ell_1(R(\lambda,B)v_1)}.
$$
And hence
$$
\mathrm{tr}(R(\lambda,A_1)-R(\lambda,B))=\ell(R(\lambda,B)v_1)=\frac{\ell_1(R(\lambda,B)^2v_1)}{1-\ell_1(R(\lambda,B)v_1)}\eqno(23)
$$
(it should be mentioned that  the condition $A_1,B\in \mathcal{N}(X)$ was not used in the proof of formula (23)).

Next, for $n\geq N$ and  $\lambda\in \Omega_{A,B},$ $|\lambda|>\lambda_0$ we have
$$
\mathrm{tr}(R(\lambda,A_n)-R(\lambda,B))=\sum\limits_{k=1}^n\mathrm{tr}(R(\lambda,A_k)-R(\lambda,A_{k-1})),
$$
where $A_k-A_{k-1}=\ell_k\otimes v_k,$ one-dimensional operator.  Now formula (23) yields
$$
\mathrm{tr}(R(\lambda,A_n)-R(\lambda,B))=
\sum\limits_{k=1}^n\frac{\ell_k(R(\lambda,A_{k-1})^2v_k)}{1-\ell_k(R(\lambda,A_{k-1})v_k)}.
$$
If we put $d_k(\lambda):=1-\ell_k(R(\lambda,A_{k-1})v_k),$  $D_n(\lambda):=\prod_{k=1}^nd_k(\lambda)$  ($\lambda\in\Omega_{A,B},$ $|\lambda|>\lambda_0$) the last formula takes the form
$$
\mathrm{tr}(R(\lambda,A_n)-R(\lambda,B))=\sum\limits_{k=1}^n\frac{d'_k(\lambda)}{d_k(\lambda)}=\frac{D_n'(\lambda)}{D_n(\lambda)}.
$$
Comparing this with corollary 4.8, we get $\Delta'_{B/A_n}/\Delta_{B/A_n}=D_n'/D_n$ and therefore $\Delta_{B/A_n}(\lambda)=C_nD_n(\lambda)$ for some constant $C_n>0.$ To show that $C_n=1$ note that $\Delta_{B/A_n}(\lambda)=\exp(\xi_{A_n,B}(\lambda))\to 1$  for $\lambda\to +\infty.$ On the other hand,
$R(\lambda,B)\to 0$  for $|\lambda|\to+\infty.$ Moreover,  $\|R(\lambda,B)\|\|S_{k-1}\|< 1$ for $\lambda\in \Omega_{A,B},$ $|\lambda|$ sufficiently large  and then by \cite[Theorem IV.1.16, Remark IV.1.17]{Kato}
$$
\|R(\lambda,A_{k-1})\|\leq \frac{\|R(\lambda,B)\|}{1-\|R(\lambda,B)\|\|S_{k-1}\|}\to 0 \ (|\lambda|\to+\infty).\eqno(24)
$$
This implies that $\ell_k(R(\lambda,A_{k-1})v_k)\to 0,$ and  $D_n(\lambda)=\prod_{k=1}^nd_k(\lambda)\to 1$ as $|\lambda|\to+\infty.$ So, $C_n=1$ and hence $\Delta_{B/A_n}(\lambda)=D_n(\lambda)$ which is equivalent to formula (21).

Next,  $A_{k-1}\in \mathcal{N}(X)$ for all $k> N$ and we have
$$
|\ell_k(R(\lambda,A_{k-1})v_k|\leq \frac{M'_{A_{k-1}}}{1+|\lambda|}\|\ell_k\|\|v_k\|\leq \frac{M'_{A_{k-1}}}{1+|\lambda|}\sum\limits_{k=1}^\infty\|\ell_k\|\|v_k\|.
$$
 But lemma 4.9 implies that $M'_{A_n}= M'_A(1-M'_A\|A-A_n\|)^{-1}\to M'_A$ as  $n\to \infty$ and therefore the sequence $M'_{A_n}$ is bounded. It follows that there exists such $\lambda_1>0$ that for $k> N$ and $|\lambda|>\lambda_1$ the inequality $|\ell_k(R(\lambda,A_{k-1})v_k|<1$ holds.
This inequality holds also for all $k\leq N$ and sufficiently large $|\lambda|,$
since $\ell_k(R(\lambda,A_{k-1})v_k)\to 0$ as $|\lambda|\to+\infty.$ Thus, $1-\ell_k(R(\lambda,A_{k-1})v_k)$ lies in the right half-plane for all $k$ and for all $\lambda\in \Omega_{A,B}$ with sufficiently large $|\lambda|.$
Now it follows from (21)  that for $\lambda\in \Omega_{A,B}$ with sufficiently large $|\lambda|$  and for sufficiently large $n$
$$
\xi_{A_n,B}(\lambda)=\sum\limits_{k=1}^n\log(1-\ell_k(R(\lambda,A_{k-1})v_k))\quad \eqno(25)
$$
($\log$ denotes the branch of the logarithm in the right half-plane that satisfies $\log 1=0$).

On the other hand, $\xi_{A_n,B}(\lambda)=\mathrm{tr}(\psi_\lambda(B)-\psi_\lambda(A_n))$\ ($\lambda>0$) and hence  by corollary 3.7
$$
|\xi_{A,B}(\lambda)-\xi_{A_n,B}(\lambda)|=|\mathrm{tr}(\psi_\lambda(A)-\psi_\lambda(A_n))|\leq \|\psi_\lambda(A)-\psi_\lambda(A_n)\|_{\mathfrak{S}_1}\to 0
$$
 as  $n\to \infty$ because the sequence  $M'_{A_n}$   is bounded.
 So, $\xi_{A_n,B}(\lambda)\to \xi_{A,B}(\lambda)$ as  $n\to \infty$ and for real and
 sufficiently large $\lambda\in \Omega_{A,B}$ the result follows from (25). The general case is valid
 because of the analyticity of both parts of the equality which we prove.  The absolute convergence of the
  series follows from the inequality $|\ell_k(R(\lambda,A_{k-1})v_k|\leq M'_{A_{k-1}}\|\ell_k\|\|v_k\|$
   that holds for sufficiently large $k$  and from the boundedness of the sequence $M'_{A_k}$.

The following corollary shows that the  perturbation determinant  for the pair  $(A,B)$ of negative operators on Banach space with nuclear difference possesses the basic properties of the  classical perturbation determinant in the Hilbert space setting (see \cite{GK}, or
 \cite[Section 8.1 ]{Ya}).

 \begin{cor} Let the conditions of theorem 4.10 are fulfilled. Then

 1) $$
\Delta_{B/A}(\lambda)=\prod\limits_{k=1}^\infty(1-\ell_k(R(\lambda,A_{k-1})v_k))
$$
for $\lambda\in \Omega_{A,B}$ with sufficiently large $|\lambda|,$ and the infinite  product converges absolutely;

2) $\Delta_{B/A}$ is analytic in $\Omega_{A,B}$;

3) $\Delta_{B/A}(z)\to 1$ as $z\to\infty$ in $\Omega_{A,B}$;

4)If operators $A, B, C\in \mathcal{N}(X)$ be such that $A-B$ and $B-C$ are nuclear, then
$$
\Delta_{B/A}(z)\Delta_{C/B}(z)=\Delta_{C/A}(z),\quad z\in \Omega_{A,B}\cap\Omega_{B,C}.
$$

In particular,
$$
\Delta_{B/A}(z)\Delta_{A/B}(z)=1,\quad z\in \Omega_{A,B}.
$$

5)Let the Banach space $X$ has the property that the trace   on  $\mathfrak{S}_1(X)$ is nilpotent in a sense  that  $\mathrm{tr}(N)=0$ for every nilpotent operator $N.$ Suppose $z_1$ is a regular point or a normal eigenvalue of the operators
$B$ and $A$ of finite algebraic multiplicities $k_0$ and $k.$ Then at the point $z_1$ the
function $\Delta_{B/A}(z)$ has a pole (or zero) of order $k_0- k$ (respectively of order $k - k_0$).
\end{cor}

Proof. 1) This follows from the definition 4.6  and theorem 4.10.

Statements 2) and 3)   follow from the definition 4.6 and corresponding  properties of $\xi_{A,B}.$

4) This follows from  theorem 4.10 and analyticity of the perturbation determinant.

5) Due to the corollary 4.8 the proof of this assertion is similar to the proof of the property 4 of the perturbation determinant in \cite[p. 267]{Ya}.

For nonpositive operators we have the following

\begin{thm}
Let the Banach space $X$ has the approximation property.  Let  $A$ and $B$ from $\mathcal{NP}(X)$
be such that $D(A)\subseteq D(B)$ and $A-B$ is nuclear.
Then $\xi_{A,B}$ has an analytic continuation into the closure of some sector $S_\theta$ with $\theta\in(0,\pi/2)$ and for every negative operator monotone function $\varphi,$ $\varphi'(-0)\ne \infty,$ with the property
$$
\int\limits_0^{\infty}\frac{|\varphi(-x)|}{1+x^2}dx<\infty\eqno(\ast)
$$
the following trace formula holds:
$$
\mathrm{tr}(\varphi(A)-\varphi(B))=\lim\limits_{\varepsilon\to +0}\frac{1}{2\pi i}\int\limits_{\partial S_\theta}\xi_{A,B}(z+\varepsilon)\varphi'(z)dz.\eqno(26)
$$
\end{thm}

Proof.
By theorem 4.2 for every $\varepsilon>0$
$$
\mathrm{tr}(\varphi(A-\varepsilon I)-\varphi(B-\varepsilon I))=\frac{1}{2\pi i}\int\limits_{\Gamma_{A-\varepsilon I,B-\varepsilon I}}\xi_{A-\varepsilon I,B-\varepsilon I}(z)\varphi'(z)dz.
$$
 Note that  $M_{A-\varepsilon I}\leq 2M_A$ and $M_{B-\varepsilon I}\leq 2M_B$ since, for example,
$$
M_{A-\varepsilon I}=\sup\limits_{\lambda>0}\|\lambda((\lambda+\varepsilon)I-A)^{-1}\|\leq
 \sup\limits_{\lambda>0}\|(\lambda+\varepsilon)((\lambda+\varepsilon)I-A)^{-1}\|+
 $$
 $$
 \varepsilon\sup\limits_{\lambda>0}\|((\lambda+\varepsilon)I-A)^{-1}\|\leq M_A+\varepsilon\sup\limits_{\lambda>0}\frac{M_A}{\lambda+\varepsilon}=2M_A.
$$
Formula (7) and the second resolvent identity imply that
$$
(\varphi(A)-\varphi(B))-(\varphi(A-\varepsilon I)-\varphi(B-\varepsilon I)) =
$$
$$
  \int\limits_{(0,\infty)}(R(t,A)(A-B)R(t,B)-R(t,A-\varepsilon I)(A-B)R(t,B-\varepsilon I))td\mu(t)=
$$
$$
 \int\limits_{(0,\infty)}R(t,A)(A-B)(R(t,B)-R(t,B-\varepsilon I))td\mu(t)+
 $$
 $$
 \int\limits_{(0,\infty)}(R(t,A)-R(t,A-\varepsilon I))(A-B)R(t,B-\varepsilon I)td\mu(t)=
$$
$$
\varepsilon \int\limits_{(0,\infty)}R(t,A)(A-B)R(t,B)R(t,B-\varepsilon I)td\mu(t)+
 $$
  $$
 \varepsilon \int\limits_{(0,\infty)}R(t,A)R(t,A-\varepsilon I)(A-B)R(t,B-\varepsilon I)td\mu(t).
$$
Since $\|R(t,A-\varepsilon I)\|\leq 2M_A/(1+t),$   $\|R(t,B-\varepsilon I)\|\leq 2M_B/(1+t),$ and $\int_{(0,\infty)}d\mu(t)/t<\infty$ (see corollary 3.2)  the last equality implies in view of (ii) $$
\|(\varphi(A)-\varphi(B))-(\varphi(A-\varepsilon I)-\varphi(B-\varepsilon I))\|_{\mathfrak{S}_1}\to 0 \mbox{ as } \varepsilon\to +0.
$$
 Then
$$
\mathrm{tr}(\varphi(A)-\varphi(B))=\lim\limits_{\varepsilon\to +0}\frac{1}{2\pi i}\int\limits_{\Gamma_{A-\varepsilon I,B-\varepsilon I}}\xi_{A-\varepsilon I,B-\varepsilon I}(z)\varphi'(z)dz.
$$

Recall that  $\Omega_{A-\varepsilon I,B-\varepsilon I}$ is any set of the form  $S_\theta\cup B_\delta(0)$ with $0<\theta<\max\{\arcsin1/M_{A-\varepsilon I}, \arcsin1/M_{B-\varepsilon I}\}$ such that the condition ($\mathrm{ii}'$) from the Preliminaries holds. So, one can take $0<\theta<\arcsin1/2\max\{M_{A},M_{B}\}.$ And since the function $\xi_{A-\varepsilon I,B-\varepsilon I}$ is holomorphic in some neighborhood of the closure of $\Omega_{A-\varepsilon I,B-\varepsilon I},$ the Cauchy Theorem implies
$$
\mathrm{tr}(\varphi(A)-\varphi(B))=\lim\limits_{\varepsilon\to +0}\frac{1}{2\pi i}\int\limits_{\partial S_\theta}\xi_{A-\varepsilon I,B-\varepsilon I}(z)\varphi'(z)dz.
$$

Moreover,  the condition (ii) yields that $|\eta_{A,B}(z)|\leq C/|z|^2$ for $z\in S_\theta$  ($C=M'_AM'_B\|A-B\|_{\mathfrak{S}_1}$). Therefore formula (10) gives  the analytic continuation of $\xi_{A,B}$  into  the closure of  $S_\theta$ as in the proof of theorem 4.2. Consequently the equalities
$\eta_{A-\varepsilon I,B-\varepsilon I}(z)=\mathrm{tr}(R(z+\varepsilon,A)-R(z+\varepsilon,B)),$ and $\varepsilon+L_z=L_{z+\varepsilon}$
imply for $z\in S_\theta$
$$
\xi_{A-\varepsilon I,B-\varepsilon I}(z)=
$$
$$
-\int\limits_{L_z}\eta_{A-\varepsilon I,B-\varepsilon I}(\zeta)d\zeta=-\int\limits_{\varepsilon+L_z}\mathrm{tr}(R(\zeta,A)-R(\zeta,B))d\zeta=\xi_{A,B}(z+\varepsilon)
$$
and the result follows.

\begin{rem}
Passage to the limit under the integral sign in (26) is impossible since  the  integral resulting from this may diverge. Indeed, let $X=\mathbb{C}^2,$  $A=-I,$ and $Bx=(-x_1,0).$ Then  $\xi_{A,B}(z)=\log(1+1/z)\sim 1/z$ as $z\to 0$ but $\varphi'(-0)>0$ for every $\varphi\in\mathcal{ OM}_-, \varphi\ne 0.$
\end{rem}

\subsection{The case of affine functions}
The  formula (LK) does not valid for  affine functions  $\varphi(s)=c+bs$ as theorem 4.2 shows. In this subsection we  prove that Lifshitz-Kre\u{\i}n trace formula remains valid for affine functions  if the integral  is understood  in  some generalized sense.

\begin{thm}  Let the Banach space $X$ has the approximation property. For any operators $A,B\in \mathcal{N}(X)$ such that $D(A)\subseteq D(B)$ and
operator $A- B$ is nuclear the following equality is valid:
$$
\mathrm{tr}(A-B)=\lim\limits_{\lambda\to+\infty}\frac{1}{2\pi i}\int\limits_{\Gamma_{A,B}}\frac{\lambda^2}{(\lambda-z)^2}\xi_{A,B}(z)dz.
$$
\end{thm}

Proof. In the proof we use notation and facts from the proof of theorem 4.10. Let $A- B$ has the form (20). Formula (23) implies for $k\in \mathbb{N}, \lambda>\lambda_0$ that
$$
\lambda^2\mathrm{tr}(R(\lambda,A_k)-R(\lambda,A_{k-1}))=
\frac{\ell_1((\lambda R(\lambda,A_{k-1}))^2v_k)}{1-\ell_1(R(\lambda,A_{k-1})v_k)}\eqno(27)
$$
($R(\lambda,A_k)$ exists for $\lambda>\lambda_0$, see the proof of theorem 4.10). Moreover,
since $\lambda R(\lambda,A_{k-1})=I+A_{k-1}R(\lambda,A_{k-1}),$
we have for all $x\in X$
$$
(\lambda R(\lambda,A_{k-1}))^2x=x+2A_{k-1}R(\lambda,A_{k-1})x+(A_{k-1} R(\lambda,A_{k-1}))^2x.\eqno(28)
$$
If $x\in D(A_{k-1})$ then, by (24),
$$
A_{k-1}R(\lambda,A_{k-1})x=R(\lambda,A_{k-1})A_{k-1}x\to 0\ \mbox{ as } \lambda\to+\infty.\eqno(29)
 $$
 Since $D(A_{k-1})=D(B)$ is dense in $X,$ to prove (29) for an arbitrary  $x\in X$ it suffices to show that for every $k$ the family of bounded operators $(A_{k-1}R(\lambda,A_{k-1}))_{\lambda>\lambda_0}$ is uniformly bounded. To this end note that    we have from (24)  for $\lambda>\max\{\lambda_0,1\}$ that
 $$
 \|A_{k-1}R(\lambda,A_{k-1})\|=
 $$
 $$
 \|I+\lambda R(\lambda,A_{k-1})\|\leq 1+\lambda\|R(\lambda,A_{k-1})\|\leq 1+\frac{\lambda\|R(\lambda,B)\|}{1-\|V_{k-1}\|\|R(\lambda,B)\|}.
 $$
Since $\|V_{k-1}\|\|R(\lambda,B)\|\leq M_B\|V_{k-1}\|/(1+\lambda)\leq\lambda_0/(1+\lambda),$ it follows that
$$
 \|A_{k-1}R(\lambda,A_{k-1})\|\leq 1+\frac{\lambda M_B}{1+\lambda-\lambda_0}\leq 1+\lambda_0 M_B.
 $$
In turn, it follows that for every $x\in X$
$$
(A_{k-1}R(\lambda,A_{k-1}))^2x=A_{k-1}R(\lambda,A_{k-1})A_{k-1}R(\lambda,A_{k-1})x\to 0 \mbox{ as } \lambda\to+\infty.
$$
Therefore taking into account  (28) we have ($S_0:=O$)
$$
\lim\limits_{\lambda\to +\infty}\ell_k((\lambda R(\lambda,A_{k-1}))^2v_k)=
$$
$$
\lim\limits_{\lambda\to +\infty}\ell_k(v_k)+2\ell_k(A_{k-1}R(\lambda,A_{k-1})v_k)+\ell_k((A_{k-1}R(\lambda,A_{k-1}))^2v_k)=
$$
$$
\ell_k(v_k)=\mathrm{tr}(V_k-V_{k-1}).
$$
Now (27) yields that for all $k=1,\dots,n$
$$
\lim\limits_{\lambda\to+\infty}\lambda^2\mathrm{tr}(R(\lambda,A_k)-R(\lambda,A_{k-1}))=\mathrm{tr}(S_k-S_{k-1}).
$$
Summing this equations  we get
$$
\lim\limits_{\lambda\to+\infty}\lambda^2\mathrm{tr}(R(\lambda,A_n)-R(\lambda,B))=\mathrm{tr}S_n=\mathrm{tr}(A_n-B).
$$
Since $\lim\limits_{n\to \infty}\mathrm{tr}(A_n-B)=\mathrm{tr}(A-B),$  it follows that
$$
\mathrm{tr}(A-B)=\lim\limits_{n\to \infty}\lim\limits_{\lambda\to+\infty}\lambda^2\mathrm{tr}(R(\lambda,A_n)-R(\lambda,B)).
$$
On the other hand,
$$
\lambda^2\mathrm{tr}(R(\lambda,A_n)-R(\lambda,B))\to \lambda^2\mathrm{tr}(R(\lambda,A)-R(\lambda,B)) \mbox{ as } n\to\infty
$$
uniformly with respect to $\lambda.$ Indeed,
$$
\lambda^2|\mathrm{tr}(R(\lambda,A_n)-R(\lambda,B))- \mathrm{tr}(R(\lambda,A)-R(\lambda,B))|=
$$
$$
\lambda^2|\mathrm{tr}(R(\lambda,A_n)-R(\lambda,A))|\leq
\lambda^2\|R(\lambda,A_n)-R(\lambda,A)\|_{\mathfrak{S}_1}\leq
$$
$$
\lambda^2\|R(\lambda,A_n)\|\|R(\lambda,A)\|\|A_n-A\|_{\mathfrak{S}_1}\leq
 \frac{\lambda^2}{(1+\lambda)^2}M_{A_n}M_A\|A_n-A\|_{\mathfrak{S}_1}\leq
 $$
 $$
 \mathrm{const}\|A_n-A\|_{\mathfrak{S}_1}\to 0 \mbox{ as } n\to\infty,
$$
since the sequence $M_{A_n}$ is bounded as shown in the proof of theorem 4.10.

Hence in view of corollary 4.5 we get
$$
\mathrm{tr}(A-B)=\lim\limits_{\lambda\to+\infty}\lambda^2\lim\limits_{n\to\infty}\mathrm{tr}(R(\lambda,A_n)-R(\lambda,B))=
$$
$$
\lim\limits_{\lambda\to+\infty}\lambda^2\mathrm{tr}(R(\lambda,A)-R(\lambda,B))=\lim\limits_{\lambda\to+\infty}\frac{1}{2\pi i}\int\limits_{\Gamma_{A,B}}\frac{\lambda^2}{(\lambda-z)^2}\xi_{A,B}(z)dz
$$
as was to be proven.

\begin{cor}
Let the Banach space $X$ has the approximation property. For any operators $A,B\in \mathcal{NP}(X)$ such that $D(A)\subseteq D(B)$ and
operator $A - B$ is nuclear there is some $\theta\in (0,\pi/2)$ such that for every $\varepsilon>0$ the following equality is valid:
$$
\mathrm{tr}(A-B)=\lim\limits_{\lambda\to+\infty}\frac{1}{2\pi i}\int\limits_{\partial S_\theta}\frac{\lambda^2}{(\lambda-z)^2}\xi_{A,B}(z+\varepsilon)dz.
$$
\end{cor}

Proof. Since $A-\varepsilon I,B-\varepsilon I\in \mathcal{N}(X)$ for every $\varepsilon>0,$ we have by theorem 4.14
$$
\mathrm{tr}(A-B)=\mathrm{tr}((A-\varepsilon I)-(B-\varepsilon I))=
$$
$$
\lim\limits_{\lambda\to+\infty}\frac{1}{2\pi i}\int\limits_{\Gamma_{A-\varepsilon I,B-\varepsilon I}}\frac{\lambda^2}{(\lambda-z)^2}\xi_{A-\varepsilon I,B-\varepsilon I}(z)dz.
$$
 If we  take $0<\theta<\arcsin1/2\max\{M_{A},M_{B}\},$ then by the Cauchy Theorem  we get as in the proof of theorem 4.12
$$
\mathrm{tr}(A-B)=\lim\limits_{\lambda\to+\infty}\frac{1}{2\pi i}\int\limits_{\partial S_\theta}\frac{\lambda^2}{(\lambda-z)^2}\xi_{A-\varepsilon I,B-\varepsilon I}(z)dz=
$$
$$
\lim\limits_{\lambda\to+\infty}\frac{1}{2\pi i}
\int\limits_{\partial S_\theta}\frac{\lambda^2}{(\lambda-z)^2}\xi_{A,B}(z+\varepsilon)dz.
$$

\textbf{Acknowledgments.}

Published electronically in: Complex Analysis and Operator Theory, https://doi.org/10.1007/s11785-019-00902-5.

\end{document}